\documentclass[twoside]{article}
\usepackage{amsthm,amssymb,amscd}
\theoremstyle{plain}
 
\newtheorem{lem}{Lemma}[section]

\newtheorem{thm}{Theorem}[section]
\newtheorem{defi}{Definition}[section]
\newcommand{\bb}[1]{\mbox{$\mathbb{#1}$}}
\theoremstyle{definition}
\newtheorem{ex}{Example}[section]
\theoremstyle{remark}
\newtheorem{rem}{Remark}[section]

\title{A general construction of partial Grothendieck transformations}
\author{J\"{o}rg Sch\"{u}rmann}
\date{}

\begin{document}

\maketitle
\bibliographystyle{plain}

$ $
$ $

\begin{abstract}
Fulton and MacPherson introduced
the notion of {\em bivariant theories} related to
{\em Riemann-Roch-theorems}, especially in the context of singular spaces.
This is powerful formalism, which is a simultaneous generalization of a pair of contravariant
and covariant theories. Natural transformations of bivariant theories are called
{\em Grothendieck transformations}, and these generalize a pair of ordinary natural transformations.

But there are many situations, where such a bivariant theory or a corresponding 
Grothendieck transformation is only "partially known":
characteristic classes of singular spaces (e.g. Stiefel-Whitney or Chern classes), 
cohomology operations (e.g. singular Adams Riemann-Roch and Steenrod operations
for Chow groups) or equivariant theories (e.g. 
Lefschetz Riemann-Roch). We introduce in this paper a
simpler notion of {\em partial (weak) bivariant theories} and {\em partial Grothendieck transformations},
which applies to all these examples. Our main theorem shows, that a natural transformation of
covariant theories, which commutes with exterior products, automatically {\em extends uniquely}
to such a
partial Grothendieck transformations of suitable partial (weak) bivariant theories !
In the above geometric situations one has for example to consider only morphisms, whose target is a
{\em smooth} manifold, or more generally, a suitable "{\em homology} manifold" 
(and in the general bivariant language this is related to the existence of suitable
{\em strong orientations}). We illustrate our main theorem for the examples above, relating it to
corresponding known Riemann-Roch theorems.
\end{abstract}

\subsection*{Acknowledgements.}
This paper should be seen as a continuation of the basic work of W.Fulton
and R.MacPherson \cite{FM} about {\em bivariant theories and Grothendieck transformations}.
Its origin goes back to the work \cite{Br} of J.P.Brasselet and \cite{Y1}... \cite{Y6} of S.Yokura
on the theory of {\em bivariant Chern classes}. In particular, our main theorem goes back to a
result of S.Yokura  \cite{Y6} about the construction (and uniqueness) of a suitable (in our notion "partial")
bivariant theory of such Chern classes. The author realized the abstract bivariant background
of this construction, and our aim is to show by some important {\em examples} the power of this
abstract {\em bivariant} result! \\
In the examples of this paper we consider only the case of
{\em partial} Grothendieck transformations between bivariant theories in the sense of Fulton-MacPherson.
The general construction of new {\em partial bivariant theories} will be explained elsewhere.
The author would like to thank J.P.Brasselet and S.Yokura for some
remarks on this work. 

\section{Uniqueness of bivariant transformations}

First we recall some notions of Fulton-MacPherson \cite{FM} for the general theory
of bivariant theories and Grothendieck transformations.\\  

Let $C$ be a category with a final object $pt$ and fiber-products.
Fix in addition a class of {\em confined} maps (closed under composition and base-change, containing
all identity maps), and a class of {\em independent squares} (which we always assume to be 
fiber squares), closed under "vertical" and "horizontal" composition (as in \cite[p.17]{FM}) and containing
any square (or its "transpose") of the form
\begin{displaymath} \begin{CD}
X @>id >> X \\
@V f VV  @VV f V \\
Y @> id >> Y \:. 
\end{CD} \end{displaymath}

A bivariant theory $\bb{B}$ on the category $C$ 
assigns to each morphism
\[f: X\to Y\]
in $C$ an abelian group 
\[\bb{B}(f:X\to Y)\]
together with three (linear) operations:
\begin{enumerate}
\item {\bf product:} For morphisms $f: X\to Y$ and $g: Y\to Z$ a product
\[\bullet: \bb{B}(f: X\to Y)\times \bb{B}(g: Y\to Z) \to \bb{B}(g\circ f: X\to Z)\:.\]
\item {\bf push-down:} For morphisms $f: X\to Y$ and $g: Y\to Z$ with $f$ {\em confined}
a push-forward 
\[f_{*}: \bb{B}(g\circ f: X\to Z) \to \bb{B}(g: Y\to Z)\:.\]
\item {\bf pull-back:} For any {\em independent square} 
\begin{displaymath} \begin{CD}
X' @>g' >> X \\
@V f' VV  @VV f V \\
Y' @> g >> Y 
\end{CD} \end{displaymath}
a pull-back
\[g^{*}: \bb{B}(f: X\to Y) \to \bb{B}(f': X'\to Y')\:.\]
\end{enumerate}
These three operations are required to satisfy the following seven compabilities
(compare \cite[PartI, 2]{FM} for details):
\begin{itemize}
\item {\bf (A1)} Product is associative.
\item {\bf (A2)} Push-down is functorial.
\item {\bf (A3)} Pull-back is functorial.
\item {\bf (A12)} Product and push-down commute.
\item {\bf (A13)} Product and pull-back commute.
\item {\bf (A23)} Push-down and pull-back commute.
\item {\bf (A123)} The (bivariant) projection formula.
\end{itemize}

The bivariant theory $\bb{B}$ is by definition {\em graded} ($\bb{Z}_{2}-$ or $\bb{Z}$-graded), if 
each group $\bb{B}(f: X\to Y)$ has a grading 
\[\bb{B}^{i}(f: X\to Y)\:,\]
which is stable under push-down and pull-back as above and additive under the (bivariant) product:
\[\bullet: \bb{B}^{i}(f: X\to Y)\times \bb{B}^{j}(g: Y\to Z) \to \bb{B}^{i+j}(g\circ f: X\to Z)\:.\]

As in \cite{FM}, we use for $\alpha\in \bb{B}(f: X\to Y)$ in bivariant diagrams
a symbol near the arrow of the morphism:
\begin{displaymath} \begin{CD}
X @> \alpha > f > Y\:.
\end{CD} \end{displaymath}

We assume that the bivariant theory $\bb{B}$  has a {\em unit}
\cite[p.22]{FM} (i.e. an element $1_{X}\in \bb{B}(id:X\to X)$, which is a unit with respect to
all possible bivariant products, with $g^{*}(1_{X})=1_{X'}$ for any map $g: X'\to X$).\\

Finally, a (graded) bivariant theory $\bb{B}$ is called {\em (skew-)commutative} if,
whenever the "transpose" of an independent square 
\begin{displaymath} \begin{CD}
X' @>g' >> X \\
@V f' VV  @V f V \alpha V \\
Y' @> g >\beta > Y 
\end{CD} \end{displaymath}
is also independent, then
\[g^{*}(\alpha)\bullet \beta = f^{*}(\beta)\bullet \alpha\:\]
(or 
\[g^{*}(\alpha)\bullet \beta = (-1)^{|\alpha|\cdot |\beta|} f^{*}(\beta)\bullet \alpha\:\]
in the skew-commutative case. Here $|\cdot|$ denotes the degree).\\

Moreover, $\bb{B}^{*}(X):=\bb{B}(id: X\to X)$ is by definition the {\em associated contravariant
theory}, with a {\em cup-product} 
\[\cup:  \bb{B}^{*}(X)\times \bb{B}^{*}(X) \to \bb{B}^{*}(X)\]
induced form the bivariant product for the
composition $X\stackrel{id}{\to}X\stackrel{id}{\to}X$ (with possible grading 
$\bb{B}^{i}(X):=\bb{B}^{i}(id: X\to X)$). The {\em associated covariant theory} is defined as
$\bb{B}_{*}(X):=\bb{B}(X\to pt)$ (which is covariant for confined maps), 
with a {\em cap-product} 
\[\cap:  \bb{B}^{*}(X)\times \bb{B}_{*}(X) \to \bb{B}_{*}(X)\]
 induced form the bivariant product for the
composition $X\stackrel{id}{\to}X \to pt$ (with possible grading 
$\bb{B}_{i}(X):=\bb{B}^{-i}(X\to pt)$).\\

Suppose $C$ and $\bar{C}$ are categories with classes of confined maps, independent squares
and a final object, and consider a functor
\[\bar{ }: C \to \bar{C} \]
respecting these structures. Then a {\em Grothendieck transformation}
$\gamma$ of bivariant theories $F$ on $C$ and $H$ on $\bar{C}$
is a collection of homomorphisms
\[\gamma_{f}: F(f: X\to Y) \to H(\bar{f}: \bar{X}\to \bar{Y}) \:,\]
which  commutes with product, push-down
and pull-back operations.

\begin{rem} Note that in general a Grothendieck transformation $\gamma$ need {\em not}
preserve a possible grading of the bivariant theories!
\end{rem} 

Let us now discuss the uniqueness problem of a Grothendieck transformation.
Since $\gamma$ respects the bivariant product, one gets
\begin{equation} \label{eq:product}
\gamma_{g\circ f}(\alpha\bullet \beta) = \gamma_{f}(\alpha)\,\bar{\bullet}\, \gamma_{g}(\beta)
\end{equation}
for all $\alpha\in F(f:X\to Y)$ and $\beta\in F(g:Y\to Z)$.
So if we have a distinguished element $e_{g}\in F(g:Y\to Z)$ such that
$\gamma_{g}(e_{g})$ is a {\em strong orientation} for $\bar{g}$, i.e.
\[\,\bar{\bullet}\, \gamma_{g}(e_{g}): H(f':X'\to \bar{Y}) \to H(\bar{g}\circ f':X'\to \bar{Z})\]
is an isomorphism for all morphisms $f':X'\to \bar{Y}$ in $\bar{C}$
(compare \cite[Def.2.7.1, p.29]{FM}), then $\gamma_{f}(\alpha)$ is uniquely determined
by $\gamma_{g\circ f}(\alpha\bullet e_{g})$,
$\gamma_{g}(e_{g})$ and equation (\ref{eq:product}).\\

In particular, if the corresponding covariant transformation 
\[\gamma_{*}: F_{*}\to H_{*}\]
of the associated {\em covariant} functors is unique, then the {\em bivariant} transformation
\[\gamma_{f}: F(f:X\to Y) \to H(\bar{f}:\bar{X}\to \bar{Y})\]
is also unique provided there is some element 
$e_{Y}\in F_{*}(Y):=F(Y\to pt)$ such that
$\gamma_{*}(e_{Y})$ is a {\em strong orientation}.\\
 
The important example for us comes from the theory of characteristic classes of singular spaces,
with $Y$ a {\em smooth} manifold. Here the confined maps are the {\em proper} maps
(as in the following examples).
By resolution of singularities, the corresponding covariant transformation 
$\gamma_{*}: F_{*}\to H_{*}$ is unique, and one has the normalization 
\begin{equation} \label{eq:norm}
\gamma_{*}(e_{Y}) = c^{*}(\bar{Y})\cap [\bar{Y}]\;,
\end{equation}
with $c^{*}(\bar{Y}):=c^{*}(T\bar{Y})\in H^{*}(\bar{Y}):=H(id:\bar{Y}\to \bar{Y})$ the
 corresponding {\em characteristic class} of the tangent bundle
$T\bar{Y}$, and $[\bar{Y}]\in H_{*}(\bar{Y}):=H(\bar{Y}\to pt)$ the {\em fundamental class} of the 
manifold $\bar{Y}$.\\

Here one can work in one of the following cases:
\begin{ex} \label{ex:charclass} In the first examples $F$ and $H$ are defined on the same
category (so that $\:\bar{ }\:$ is the identity transformation).
\begin{enumerate}
\item In the algebraic context with $F_{*}(Y):=K_{0}(Y)$ the Grothendieck group of coherent sheaves,
$e_{Y}$ the class of the structure sheaf, 
\[H_{*}(Y):=A_{*}(Y)\otimes \bb{Q}\] 
the Chow-group with
rational coefficients (resp. $H$ the operational bivariant Chow group with rational  coefficients)
and $\gamma_{*}=\tau_{*}$ the {\em Baum-Fulton-MacPherson transformation}
(compare \cite[Part II]{FM} and \cite[Chapter 17,18]{Ful}) so that $c^{*}(TY)$ is the 
total {\em Todd class} of $TY$.
\item In the algebraic context (over a ground field of characteristic zero)
with $F_{*}(Y)$ the group of algebraically  constructible functions, $e_{Y}:=1_{Y}$, 
\[H_{*}(Y):=A_{*}(Y)\] 
the Chow-group (resp. $H$ the operational bivariant Chow group)
and $\gamma_{*}=c_{*}$ the {\em Chern-Schwartz-MacPherson transformation} (\cite[Ex. 19.1.7]{Ful}
and \cite{Ken,EY1,EY2})
so that $c^{*}(TY)$ is the total {\em Chern class} of $TY$.
\item In the complex analytic (algebraic) context with $F_{*}(Y)$ the group of 
analytically (algebraically) constructible functions, $e_{Y}:=1_{Y}$, 
\[H_{2*}(Y):=H^{BM}_{2*}(Y,\bb{Z})\] 
the Borel-Moore Homology-group (resp. $H$ the bivariant Homology) in even degrees
and $\gamma_{*}=c_{*}$ the {\em Chern-Schwartz-MacPherson transformation} \cite{M,Br,BrS}
so that $c^{*}(TY)$ is the total {\em Chern class} of $TY$. If we consider only spaces which can be embedded
into complex analytic (algebraic) manifolds, then we can use instead of $F_{*}(Y)$ the isomorphic
theory $L(X)$ of (conic) Lagrangian cycles \cite{Sab1,Sab2,Ken,Sch}, with $e_{Y}$ the zero-section of 
the cotangent bundle $T^{*}Y$ for $Y$ smooth. Then $\gamma_{*}$ corresponds to the {\em Chern-Mather 
transformation} so that $c^{*}(TY)$ is again the total {\em Chern class} of $TY$.
\item In the subanalytic (or pl-) context with $F_{*}(Y)$ the group of 
subanalytically (or pl-) constructible Euler-functions with values in $\bb{Z}_{2}$, $e_{Y}:=1_{Y}$, 
\[H_{*}(Y):=H^{BM}_{*}(Y,\bb{Z}_{2})\] 
the Borel-Moore Homology-group (resp. $H$ the bivariant Homology)
and $\gamma_{*}=w_{*}$ the {\em Stiefel-Whitney transformation} (\cite{FuMC,Sch} or \cite{Su,FM,EH})
so that $c^{*}(TY):=w^{*}(TY)$ is the total {\em Stiefel-Whitney class} of $TY$.
\item In this last example we work on the category of complex spaces with a {\em real structure},
i.e. with an antiholomorphic involution
(or in the algebraic context over the ground field $\bb{R}$ so that the associated
complex spaces have an induced real structure). In this case the transformation
$\:\bar{ }\:$ is given on the objects by 
\[X \mapsto X(\bb{R})\;\]
with $X(\bb{R})$ the fixed point set of the real structure (of the associated
complex space), and on morphisms it is just restriction to these subsets
(of the associated holomorphic map). Here we consider as $F_{*}(Y)$ the group of 
analytically (algebraically) constructible functions, which are {\em invariant} 
under the real structure, $e_{Y}:=1_{Y}$ and 
\[H_{*}(\bar{Y}):= H^{BM}_{*}(Y(\bb{R}),\bb{Z}_{2})\] 
the Borel-Moore Homology-group (resp. $H$ the bivariant Homology).
The transformation $\gamma_{*}$ is given as the composition
\begin{displaymath} \begin{CD}
F(X) @> |\,mod\,2 >> F(X(\bb{R})) @> w_{*} >> H_{*}(X(\bb{R})) \,,
\end{CD} \end{displaymath}
with $|\,mod\,2$ given by restriction and taking the values $mod\,2$, and $w_{*}$ is the
{\em Stiefel-Whitney transformation} of the example before (note that $|\,mod\,2$ maps to
real analytically (algebraically) constructible Euler-functions,
and it is functorial with respect to push-down for proper equivariant holomorphic maps).
So for a smooth manifold $Y$ one has 
\[\gamma_{*}(e_{Y}) = w^{*}(TY(\bb{R}))\cap [Y(\bb{R})]\;.\quad \Box\]
\end{enumerate}
\end{ex}

Note that in all cases the {\em fundamental class} $[\bar{Y}]\in H_{*}(\bar{Y})$ is a strong orientation
and $c^{*}(T\bar{Y})\in H^{*}(\bar{Y})$ is {\em invertible} so that 
$\gamma_{*}(e_{Y}) = c^{*}(T\bar{Y})\cap [\bar{Y}]$ is also
a {\em strong orientation} ! 

\begin{rem} The above reasoning for the {\em uniqueness} of these bivariant characteristic
classes is just the argument of \cite{Y1} written up in the bivariant language.
In that paper Yokura shows in particular the uniqueness of {\em bivariant Chern classes} for morphisms
with smooth target (extending a result of \cite{Zhou,Zhou2}, comparing the bivariant Chern classes of
\cite{Br} and \cite{Sab2} for the case of a smooth one-dimensional target).
Moreover, our bivariant argument goes already back to \cite[p.70]{FM}, where they prove the 
uniqueness of the bivariant {\em Stiefel-Whitney classes} in the pl-context.
Their proof consists of three steps:
\begin{enumerate} 
\item By resolution of singularities, the corresponding covariant transformation 
$w_{*}: F_{*}\to H_{*}$ is unique \cite[end of p.70]{FM}.
\item  Uniqueness in the case of a {\em smooth} target \cite[eq. (***), p.70]{FM}, with the same argument
as we explained above.
\item  Uniqueness in the general case by reducing it to (a pullback of) the case 2.
\end{enumerate}
Only the last step does {\em not} generalize to the other bivariant characteristic classes!
Moreover, steps 1. and 2. also apply to suitable "equivariant versions" of the cases given in
example \ref{ex:charclass} (as will be explained somewhere else and compare with
\cite{EG2}).
\end{rem}  
 
\begin{rem} \label{rem:uniqsing}
Our {\em uniqueness} argument applies also to a suitable
{\em singular} space $\bar{Y}$, if for example the
fundamental class $[\bar{Y}]\in H_{*}(\bar{Y})$ is a strong orientation
with $c^{*}(\bar{Y})\in H^{*}(\bar{Y})$ {\em invertible},
where this cohomology class is {\em defined} by
\[\gamma_{*}(e_{Y}) = c^{*}(\bar{Y})\cap [\bar{Y}] \:.\]
In example \ref{ex:charclass} this is true in the case
\begin{enumerate}
\item for $Y$ an {\em Alexander scheme} in the sense of \cite{Vis,Ki2,Ki3}.
\item , if $Y$ has a {\em Chow-isomorphism} $\tilde{Y} \to Y$ in the sense
of \cite[def.3.5, p.296]{Ki1}, with $\tilde{Y}$ a {\em smooth} manifold
(e.g. $Y$ is a cusp curve \cite[ex.3.7(2), p.297]{Ki1}).
\item $Y$ is an oriented {\em $\bb{Z}$-homology manifold}.
\item $Y$ is a {\em $\bb{Z}_{2}$-homology manifold}. 
\item $Y(\bb{R})$ is a {\em $\bb{Z}_{2}$-homology manifold}.
\end{enumerate}
The fundamental class $[\bar{Y}]$ is in the last three cases 
a strong orientation by \cite[prop.7.3.2, p.85]{FM}.
Moreover, in all these five cases is $H^{*}(\bar{Y})$ a {\em commutative} ring with
unit $1_{\bar{Y}}$ such that $c^{*}(\bar{Y})- 1_{\bar{Y}}$ is {\em nilpotent},
i.e. $\gamma_{*}(e_{Y}) = [\bar{Y}]$ up to lower order terms!
\end{rem}

But even if the corresponding {\em covariant} transformation
$\gamma_{*}$ maybe is not unique, one can use the same argument to describe 
the {\em bivariant} transformation
\[\gamma_{f}: F(f:X\to Y) \to H(\bar{f}:\bar{X}\to \bar{Y})\]
in terms of $\gamma_{*}$. Important examples come from {\em cohomology operations}
(compare \cite[4.1.4,4.2, p.46-49]{FM}):
\begin{ex} \label{ex:cohop}
\begin{enumerate}
\item $\gamma: H\to H$ is the Grothendieck transformation on $\bb{Z}_{2}$-valued  
bivariant homology $H$ induced by the {\em Steenrod square} \cite[4.2.1, p.47]{FM}.
Then one has for a smooth manifold $Y$ with fundamental class $e_{Y}:=[Y]\in H_{*}(Y)$:
\[\gamma_{*}([Y]) =  c^{*}(TY)\cap [Y]\:,\]
with $c^{*}(TY) =  w^{*}(TY)^{-1}$ the {\em inverse Stiefel-Whitney class} of $TY$.
\item One looks for a bivariant transformation 
\[\gamma=\psi(j): K\to \bb{Z}[1/j]\otimes K\:,\] 
with $K$ a suitable bivariant algebraic $K$-theory (as in \cite[Part II]{FM}) such that
\[\gamma_{*}= \psi_{j}:   K_{*}\to \bb{Z}[1/j]\otimes K_{*}\]
is the covariant transformation described in the {\em singular Adams Rie\-mann-Roch theorem} 
\cite[p.190]{FL} ($j\in \bb{N}$). 
Then one has for a smooth manifold $Y$, with the fundamental class $e_{Y}:=[Y]\in K_{*}(Y)$
given by the class of the structure sheaf:
\[\gamma_{*}([Y]) =  c^{*}(TY)\cap [Y]\:,\]
with $c^{*}(TY) =  \theta^{j}(TY^{\vee})^{-1}$ the {\em inverse cannibalistic class} of the dual of $TY$.
\item We work in the algebraic context over a base field, whose characteristic is different
from a {\em fixed} prime $p$. Then one can ask for a bivariant transformation 
\[\gamma=C(p): A^{*}\otimes \bb{Z}_{p} \to A^{*}\otimes \bb{Z}_{p}\:,\]
with $A^{*}$ the bivariant Chow groups, whose associated covariant transformation
\[\gamma_{*}=C_{\bullet}: A_{*}\otimes \bb{Z}_{p} \to A_{*}\otimes \bb{Z}_{p}\]
is the total {\em Steenrod $p$-th power operation} on Chow groups recently constructed
in \cite{Bro}. By \cite[prop.9.4(iii)]{Bro} we have for $Y$ {\em smooth}:
\[S_{\bullet}([Y]) = c^{*}(TY) \cap [Y] \:,\]
with $c^{*}(TY) =  w^{*}(TY)^{-1}$ the inverse of the characteristic class
\[w^{*}(TY):= \prod \: \Bigl(1+ \lambda_{i}^{p-1}\Bigr) \]
defined in terms of the {\em Chern roots} $\lambda_{i}$ of $TY$ \cite[p.1891]{Bro}.
$ \quad \Box$
\end{enumerate}
\end{ex}

Another class of important examples comes from {\em equivariant} theories and corresponding
{\em Lefschetz Riemann-Roch theorems}:
\begin{ex} \label{ex:LefRR}
Here we consider for simplicity the complex algebraic context with schemes of finite type
over $spec(\bb{C})$, together with  an {\em automorphism of finite order} (of course the morphisms
have to commute with these automorphisms. For another treatment of the first example in the algebraic
context over an algebraically closed field compare with \cite{BFQ}).   
\begin{enumerate}
\item We consider only quasi-projective schemes and work in the context of "coherent sheaves"
given in \cite[10.1, p.103-105]{FM} with a bivariant transformation
\[\tau: K^{eq}_{alg}(f:X\to Y)\to H(|f|:|X|\to |Y|)\otimes \bb{C}\:.\] 
Here $K^{eq}_{alg}$ is the Grothendieck group of "equivariant $f$-perfect complexes". Moreover,
$\:\bar{ } \;=:\:|\:\:|\:$ is given by {\em restriction to the fixed point set} of the associated
complex spaces, and these fixed point sets are assumed to be {\em projective} !
Finally $H$ is the usual bivariant Homology.
Then one has for a smooth manifold $Y$ and $e_{Y}\in K^{eq}_{alg}(Y\to pt)$ given by the structure
sheaf (with its canonical isomorphism $\phi$ lifting the automorphism of $Y$):
\[\gamma_{*}(e_{Y}) = \Bigl(ct(\lambda_{|Y|} Y)^{-1} \cup td^{*}(T|Y|)\Bigr)\cap [|Y|]\:.\]
Here
\[\lambda_{|Y|} Y:= \sum (-1)^{i} [\Lambda^{i} N] \in K^{eq}_{alg}(id: |Y|\to |Y|)\]
is the corresponding (invertible) "Euler class" of the equivariant conormal sheaf $N$ to $|Y|$ in $Y$
(and for the definition of the "Chern trace" $ct$ compare with \cite[p.104]{FM}).
\item We work in the context of "algebraically constructible sheaves" of vector spaces (over a field $k$),
whose stalks are finite dimensional.
Consider the Grothendieck group $K^{eq}_{c}(X\to pt)$ of equivariant
algebraically constructible sheaves on $X$ and define the transformation
$\gamma_{*}: K^{eq}_{c}(X\to pt)\to H_{*}(|X|)\otimes k$ as the composition
\begin{displaymath} \begin{CD}
K^{eq}_{c}(X\to pt) @> tr_{|X|}\circ \,|\,>> F(|X|)\otimes k @> c_{*}\otimes k >>  H_{*}(|X|)\otimes k\,.
\end{CD} \end{displaymath}
Here $tr_{|X|}\circ \,|$ is given by taking stalkwise the trace of the restriction to the fixed
point set, and $c_{*}$ is the Chern class transformation on the group $F(|X|)$ of algebraically
constructible functions on $|X|$. Here we can choose for $H_{*}$ the Borel-Moore homology group
or the Chow group. One can show that $\gamma_{*}$ (or $tr_{|X|}\circ \,|$) 
is a natural transformation for equivariant
proper morphisms $f: X\to Y$, where
\[f_{*}: K^{eq}_{c}(X\to pt)\to K^{eq}_{c}(Y\to pt)\]
is defined by the (alternating sum of the classes of the equivariant) higher direct image sheaves.
So one looks for a bivariant extension of $\gamma_{*}$.
If $Y$ is a smooth manifold, one has for $e_{Y}\in K^{eq}_{c}(Y\to pt)$ given by the constant
sheaf $k_{Y}$ (with its canonical isomorphism $\phi$ lifting the automorphism of $Y$):
\[\gamma_{*}(e_{Y}) = c^{*}(T|Y|)\cap [|Y|]\,,\]
with $c^{*}(T|Y|)$ the Chern class of the tangent bundle of the fixed point manifold $|Y|$.
$ \quad \Box$
\end{enumerate}
\end{ex} 

Also we are mainly interested in applications to {\em singular spaces}, our results also apply
to suitable categories of manifolds, especially to Riemann-Roch theorems in the framework of
{\em oriented cohomology (pre)theories} as recently studied in \cite{LM,Lev,Pa,PaSm}
(and compare with \cite[p.46-52]{FM} for similar results in the context of
differentiable manifolds, which are oriented with respect to suitable cohomology theories,
e.g. for complex manifolds and {\em multiplicative complex oriented} cohomology theories):\\

An  {\em oriented cohomology (pre)theory} is a suitable (contra-variant) functor
\[A: Sm \to Rings \]
on the category $Sm$ of smooth quasi-projective varieties over a field $k$
with values in the category of (commutative graded) rings with unit.
$A$ is also covariant functorial with respect to proper morphisms
(of constant relative dimension). Note that a proper morphism of quasi-projective varieties
is projective!\\
These satisfy a projection-formula (i.e. the push-down for $f: X\to Y$ is a two-sided (!)
$A(Y)$-module operator), and  the base-change property $g^{*}f_{*}=f'_{*}g'^{*}$
for any {\em transverse} cartesian diagram
\begin{displaymath} \begin{CD}
X' @>g' >> X \\
@V f' VV  @VV f V \\
Y' @> g >> Y \:,
\end{CD} \end{displaymath}
with $f,f'$ proper (of constant relative dimension).
$A$ has in addition to satisfy some other properties (which are not important for us),
and these imply especially a corresponding theory of Chern classes with ("universally central")
values in $A$ \cite[thm.2.2.2]{Pa}.\\

As we explain later on, we therefore get on $Sm$ a {\em simple bivariant theory} $A$ defined
by $A(f:X\to Y):=A(X)$. The cofined maps are the  proper morphisms
(of constant relative dimension), and the independent squares are the 
transverse cartesian diagrams, with the obvious push-down and pull-back transformations.
Finally the bivariant product
\[\bullet: A(f:X\to Y) \times A(g:Y\to Z) \to A(g\circ f: X\to Z)\]
is just given by $\alpha\bullet \beta := \alpha  \cup f^{*}(\beta)$, with $\cup$ the given product
of the ring-structure.
Then $e_{Y}:=1_{Y}\in A(Y)=A(Y\to pt)$ is (trivially) a strong orientation.\\

If $\phi: A\to B$ is a "nice" ring morphism of two such theories (compare \cite[thm.2.5.4, p.46]{Pa}
for details), then one has a corresponding Riemann-Roch theorem saying that
the composition $\gamma_{*}$:
\begin{displaymath} \begin{CD}
A_{*}(X)=A^{*}(X) @> \phi >> B^{*}(X) @> \cup td_{\phi}(TX) >> B^{*}(X)=B_{*}(X) 
\end{CD} \end{displaymath}
is a natural transformation between the corresponding covariant theories
(i.e. commutes with push-down). Here $td_{\phi}$ is a "suitable" Todd genus associated to $\phi$
\cite[Def. 2.5.2, p.45]{Pa}. Since $\phi$ is a ring morphism we get $\phi(1_{Y})=\bar{1}_{Y}$,
and therefore
\[\gamma_{*}(e_{Y}) = c^{*}(Y)\cap [Y]\:,\]
with $c^{*}(Y):= td_{\phi}(TY)\in B(Y)=B^{*}(Y)$ and $[Y]:=\bar{1}_{Y}\in B(Y)=B_{*}(Y)$.
Moreover $c^{*}(Y):= td_{\phi}(TY)$ is by definition invertible and "universally central"
(i.e. any pullback of it is central). $\quad \Box$\\
$ $\\

Let us come back to the general bivariant context.
Then one gets in all the preceding examples the following {\em explicit} (!) description
of a corresponding {\em bivariant} transformation 
\[\gamma_{f}: F(f:X\to Y)\to H(\bar{f}:\bar{X}\to \bar{Y})\]
in the case $Y$ a smooth manifold (compare \cite[thm.3.4, thm.3.7]{Y1}):
\begin{equation} \label{eq:exbiv}
\gamma_{f}(\alpha)\,\bar{\bullet}\, [\bar{Y}] = \bar{f}^{*}(c^{*}(\bar{Y}))^{-1} \cap 
\gamma_{*}(\alpha\bullet e_{Y}) \in H_{*}(\bar{X})\;.
\end{equation}
Here $\bar{f}^{*}: H^{*}(\bar{Y})\to H^{*}(\bar{X})$ is the pullback of the associated {\em contravariant}
theory and 
\[\cap : H^{*}(\bar{X}) \times H_{*}(\bar{X}) \to H_{*}(\bar{X})\]
is the corresponding {\em cap-product} \cite[p.23]{FM}. This formula (\ref{eq:exbiv}) is also true,
with the same proof, for $Y$ a {\em singular} space as in remark \ref{rem:uniqsing} !\\

Since our bivariant theories $H$ are commutative \cite[p.22]{FM}
(or $c^{*}(Y)^{-1}$ is universally central in the context of
 oriented cohomology (pre)theories), this follows from the following 
commutative diagram (whose left square is independent):

\begin{equation} \label{eq:CD1} \begin{CD}
\bar{X}  @> \gamma_{f}(\alpha)> \bar{f}>  \bar{Y} @> [\bar{Y}] >> \bar{pt} \\
@V \bar{f}^{*}(c^{*}(\bar{Y}))^{-1}V id_{\bar{X}} V  @V id_{\bar{Y}}  Vc^{*}(\bar{Y})^{-1} V  
@VV \bar{1}_{\bar{pt}}=[\bar{pt}] V\\
\bar{X} @>\bar{f} > \gamma_{f}(\alpha)> \bar{Y} @>> \gamma_{*}(e_{Y})> \bar{pt} \;,
\end{CD} \end{equation}

since the {\em associativity} of the bivariant product implies
\[\gamma_{f}(\alpha)\,\bar{\bullet}\, [\bar{Y}] = \gamma_{f}(\alpha)\,\bar{\bullet}\, 
\Bigl( c^{*}(\bar{Y})^{-1}\,\bar{\bullet}\, \gamma_{*}(e_{Y}) \Bigr) = \]
\[ \Bigl( \gamma_{f}(\alpha)\,\bar{\bullet}\,  c^{*}(\bar{Y})^{-1} \Bigr) \,\bar{\bullet}\, \gamma_{*}(e_{Y}) =
\Bigl( \bar{f}^{*}(c^{*}(\bar{Y}))^{-1}\,\bar{\bullet}\,  \gamma_{f}(\alpha) \Bigr) 
\,\bar{\bullet}\, \gamma_{*}(e_{Y})   =\]
\[\bar{f}^{*}(c^{*}(\bar{Y}))^{-1}\,\bar{\bullet}\, \Bigl( \gamma_{f}(\alpha)  
\,\bar{\bullet}\, \gamma_{*}(e_{Y}) \Bigr)  =
\bar{f}^{*}(c^{*}(\bar{Y}))^{-1} \cap \gamma_{*}(\alpha\bullet e_{Y}) \;.\]

\begin{rem} \label{rem:verdier}
Again this is just the argument of \cite{Y1} written up in the bivariant language.
In the special case $\alpha=1_{f}:=1_{X}$ in the pl-context for the bivariant Stiefel-Whitney class,
the formula (\ref{eq:exbiv}) is already explained (implicitely) in \cite[p.12/13]{FM}:
\[w_{*}(X) = f^{*}(w^{*}(TY)) \cap \xi \;,\]
with $\xi := w_{f}(1_{f})\,\bar{\bullet}\, [Y] \in H_{*}(X,\bb{Z}_{2})$.
Similarly, if $f: X\to Y$ is {\em smooth} (so that $X$ is also smooth) one has for $e_{f}\in F(f: X\to Y)$
equal to $1_{X}$ (or the class of the structure  sheaf of $X$, resp. the constant sheaf $k_{X}$)
 the relation 
$\: e_{f}\bullet e_{Y} = e_{X}\;.$
Therefore (\ref{eq:exbiv}) implies
\[\gamma_{f}(e_{f})\,\bar{\bullet}\, [\bar{Y}] = 
f^{*}(c^{*}(T\bar{Y}))^{-1} \cap \Bigl( c^{*}(T\bar{X})) \cap [\bar{X}] \Bigr) =\]
\[\Bigl( \bar{f}^{*}(c^{*}(T\bar{Y}))^{-1} \cup  c^{*}(T\bar{X}) \Bigr) \cap [\bar{X}] = 
c^{*}(T_{\bar{f}})  \cap [\bar{X}] =\]
\[c^{*}(T_{\bar{f}}) \,\bar{\bullet}\, \Bigl( [\bar{f}]\,\bar{\bullet}\, [\bar{Y}] \Bigr) = 
\Bigl( c^{*}(T_{\bar{f}}) \,\bar{\bullet}\,  [\bar{f}] \Bigr)\,\bar{\bullet}\, [\bar{Y}] \;.\]
Here $T_{\bar{f}}$ is the {\em relative tangent bundle} and 
$[\bar{f}]\in H(\bar{f}: \bar{X}\to \bar{Y})$ the {\em relative orientation class}.
Since $[\bar{Y}]$ is a strong orientation, one gets the well known formula
\begin{equation} \label{eq:relVerdier}
\gamma_{f}(e_{f}) = c^{*}(T_{\bar{f}}) \,\bar{\bullet}\,  [\bar{f}] \;.
\end{equation}
Compare \cite[prop. 6A, p.65]{FM} for the pl-Stiefel-Whitney class,
\cite[prop.3.7]{EOY} for the Chern class and \cite[formula (*), p.124]{FM} for the Todd class
in the complex algebraic context (which is stated there more generally for local complete
intersection morphisms). In the last case our argument above works for any morphism
$f: X\to Y$ of smooth spaces $X,Y$. More generally it works for a
morphism $f$ of manifolds, if we have an element
$e_{f}\in F(f:X\to Y)$ with $e_{f}\bullet e_{Y} = e_{X}$, and a  relative orientation class
$[\bar{f}]\in H(\bar{f}: \bar{X}\to \bar{Y})$ with 
$[\bar{f}] \,\bar{\bullet}\, [\bar{Y}] = [\bar{X}]\:$.
\end{rem}

A similar diagram as in (\ref{eq:CD1}) can also be used to reduce under suitable
assumptions a {\em bivariant product} to the {\em cup- and cap-product} of the associated
contra- and covariant theory. Consider a (skew-)commutative (graded)
bivariant theory $H$ and two objects $Y,Z$
such that $H_{*}(Y)$ and $H_{*}(Z)$ contain a strong orientation $[Y]$ and $[Z]$.
Then we have in particular isomorphisms
\[\bullet [Y]: H^{*}(Y)\to H_{*}(Y) \:\:,\:\: \bullet [Y]: H(f: X\to Y)\to H_{*}(X)\]
and similarly for $Z$. Then the commutative diagram (whose left square is independent)

\begin{equation} \label{eq:CD2} \begin{CD}
X  @> \alpha> f>  Y @> \beta > g> Z \\
@V f^{*}(\beta ')V id_{X} V  @V id_{Y}  V\beta ' V  
@VV [Z] V\\
X @>f > \alpha> Y @>> [Y]> pt 
\end{CD} \end{equation}
$ $\\
implies that for $\alpha\in H(f: X\to Y)$ and $\beta\in H(g: Y\to Z)$ the 
bivariant product $\alpha\bullet \beta \in H(g\circ f: X\to Z)$ corresponds
under the isomorphism $\bullet [Z]$ to $\alpha '\cap f^{*}(\beta ')\:$,
with $\alpha ':= \alpha \bullet [Y] \in H_{*}(X)\:$:
\[(\alpha\bullet \beta) \bullet [Z] = \alpha\bullet (\beta\bullet [Z]) =
 \alpha\bullet (\beta '\bullet [Y]) = 
 (\alpha\bullet \beta ')\bullet [Y] = \] 
\[ (-1)^{|\beta '|\cdot|\alpha|}(f^{*}(\beta ')\bullet \alpha )\bullet [Y] =
(-1)^{|\beta '|\cdot|\alpha|} f^{*}(\beta ')\bullet (\alpha \bullet [Y]) =\]
\[(-1)^{|\beta '|\cdot|\alpha '|} f^{*}(\beta ')\cap \alpha '
=: \alpha '\cap  f^{*}(\beta ')\:.\]
Here $|\cdot|$ denotes the degree (which in the commutative case should be set
equal to zero in the above formula), and in the last equality we use the usual
 (graded) right module structure
of a (graded) left module. Moreover, we assume in the {\em skew-commutative} case
that all strong orientations in the above calculation 
have {\em even} degrees !\\

In the example of {\em oriented cohomology (pre)theories} this gives us back our definition
of the bivariant product $\bullet$.
Another important example is given by the bivariant Homology theory (or the operational bivariant 
Chow group) with $Y$ smooth and $[Y]$ the corresponding fundamental class 
so that the isomorphism
\[\bullet [Y]: H^{*}(Y)\to H_{*}(Y) \]
is just {\em Poincar\'{e}} or {\em Alexander duality}. Then the above equality
implies for $Z:=pt$ (and $[Z]:=1_{pt}$) :
\begin{equation}  \label{eq:cup}
\alpha ' \odot \beta  := \alpha '\cap f^{*}(\beta ')  = \alpha \bullet \beta \:.
\end{equation}
So the product $\odot: H_{*}(X)\times H_{*}(Y)\to H_{*}(X)$ corresponds under the isomorphism induced by
$\bullet [Y]$ to the bivariant product 
\[\bullet : H(f: X\to Y)\times H_{*}(Y)\to H_{*}(X)\]
(compare \cite[proof of thm.3.9]{Y3}). In particular, it is associative (if defined). \\
$ $\\

\section{Partial (weak) bivariant theories} 
We now explain (following ideas of Shoji Yokura \cite{Y6}), 
how a covariant transformation $c_{*}: F_{*}\to H_{*}$
of bivariant theories can "partially" be extended to a Grothendieck transformation of
(partial) bivariant theories. Here we introduce the following notions:

\begin{defi} \label{defi:partbiv}
Let $C$ be a category with classes of confined maps, independend squares and a final
object (as in \cite[Part I,2]{FM}). 
\begin{enumerate}
\item A {\em weak} bivariant theory $T$ assigns to each morphism
$f: X\to Y$ in $C$ an abelian group $T(f: X\to Y)$ together with three operations
{\em product}, {\em push-forward} and {\em pull-back} satisfying the compabilities $(A1)-(A23)$
as in \cite[Part II]{FM}, but not necessarily the (bivariant) projection formula $(A123)$.
\item Consider in addition a class of maps in $C$, called
{\em allowable maps}, which is closed under composition such that also the composition $g\circ f$
is allowable for any confined $f$ and allowable $g$.
A {\em partial (weak) bivariant theory} $T$ assigns to each {\em allowable} morphism
$f: X\to Y$ in $C$ an abelian group $T(f: X\to Y)$ together with three operations
{\em product}, {\em push-forward} and {\em pull-back} satisfying the compabilities $(A1)-(A123)$
(or $(A1)-(A23)$ in the weak case)
as in \cite[Part II]{FM}, but this time only for all  {\em allowable} maps: e.g. the  
{\em push-forward} 
\[f_{*}: T(g\circ f: X\to Z) \to T(g: Y\to Z)\]
is only defined for $g: Y\to Z$ {\em allowable} and
$f: X\to Y$  {\em confined}, and the {\em pull-back}
\[g^{*}: T(f: X\to Y) \to T(f': X'\to Y')\]
is only defined for any independent square
\begin{displaymath} \begin{CD}
X' @>g' >> X \\
@V f' VV  @VV f V \\
Y' @>> g > Y \:,
\end{CD} \end{displaymath}
with $f,f'$ {\em allowable} (and in the bivariant projection formula (A123) of \cite[p.21/22]{FM}
one has in addition to assume $h\circ g$ is allowable). 
\item
The partial (weak) bivariant theory $T$ is {\em graded} ($\bb{Z}_{2}-$ or $\bb{Z}$-graded), if 
each group $T(f: X\to Y)$ (for $f$ allowable) has a grading 
\[T^{i}(f: X\to Y)\:,\]
which is stable under push-down and pull-back, and additive under the (bivariant) product.
One defines the {\em (skew-)commutativity} of such a theory as in the usual bivariant
context (e.g. one asks (skew) commutativity for independent squares, whose "transpose" is also
independent, with all maps allowable).
\item
If $C,\bar{C}$ are two such categories with  allowable maps, and $\:\:\bar{ }: C \to \bar{C}$
is a functor respecting these structures (especially it transforms 
 allowable maps in $C$ to allowable maps in $\bar{C}$, then a 
{\em partial Grothendieck transformation}
$t$ of partial (weak) bivariant theories $T$ and $U$ (on $C,\bar{C}$) is a collection of homomorphisms
\[t_{f}: T(f: X\to Y) \to U(\bar{f}: \bar{X}\to \bar{Y})\]
for each allowable map $f: X\to Y$ in $C$, which commutes with product, push-forward and
pull-back operations.
\end{enumerate}
\end{defi}

All of our discussions so far (and also many of the arguments of \cite[part I, §2]{FM}) extend directly
to {\em partial} (weak) bivariant theories. Here we introduced also the "weak notions", since most of 
our arguments
in this paper work without the projection formula $(A123)$.\\

Two important differences are the following:
\begin{enumerate}
\item[(a)] In general $id_{X}: X\to X$ need not be an {\em allowable} map for an object $X$ of $C$.
Especially, one does not have in general an associated {\em contravariant} (or similarly
{\em covariant}) theory. If $X\to pt$ is {\em allowable}  for all objects $X$ of $C$, then one has 
at least a corresponding {\em covariant} theory $T_{*}$ (covariant with respect to {\em confined} maps).
\item[(b)] The {\em pull-back} $g^{*}$ maybe is also defined for a morphism $g$ which is
not {\em allowable}.
\end{enumerate}

Here are some possible ways of constructing {\em partial} bivariant theories:
\begin{enumerate}
\item Any (weak) bivariant theory $F$ on $C$ "restricts" to a {\em partial} (weak) bivariant theory,
e.g. one uses as allowable maps only those morphisms, whose target belongs to a fixed
class of objects in $C$ (containing $pt$, e.g. "smooth manifolds").
\item Similarly, consider a (partial weak)
bivariant theory $F$ on $C$ and suppose that one has for each
{\em allowable} morphism $f: X\to Y$ a subgroup $F'(f: X\to Y)$ of $F(f: X\to Y)$,
which is stable under product, push-down and pull-back. This gives then a  
{\em partial} (weak)  bivariant theory $F'$, which we  call a {\em partial subtheory} of $F$.
\end{enumerate}

\begin{ex} Let $C$ be the category of complex algebraic varieties, with the {\em proper}
maps the confined maps. Define $f: X\to Y$ to be {\em allowable}, if $Y$ is a
{\em smooth} manifold, and consider the "restrictions" $A^{*}$ and $H^{2*}$ of the
bivariant Chow groups and the bivariant Homology groups (in even degrees).
Then the {\em cycle map} $cl: A_{*}\to H_{2*}$ \cite[chap.19]{Ful} of the associated
covariant theories induces by (\ref{eq:cup}) and \cite[thm.19.2, p.380]{Ful}
a {\em partial Grothendieck transformation} $\gamma$ of these {\em partial bivariant theories}:
\begin{displaymath} \begin{CD}
\gamma_{f}: A^{*}(f:X\to Y) \simeq A_{*}(X) @> cl >> H_{2*}(X) \simeq H^{2*}(f:X\to Y) \:.
\end{CD} \end{displaymath}
Note, that a corresponding Grothendieck transformation of the origional
bivariant theories is {\em not} known.
\end{ex}

We explain now a general way how a natural transformation $c_{*}: F_{*}\to H_{*}$ of
associated {\em covariant} theories can be extended to a {\em partial Grothendieck transformation}.\\

Consider two {\em partial (weak) bivariant theories}  $F$ on $C$ and $H$ on $\bar{C}$, together with a
functor  $\:\bar{ }: C \to \bar{C}$ respecting the underlying structures.
We assume that all maps $X\to pt$ are {\em allowable} in $C$ (and the same for $\bar{C}$),
and $F_{*}(pt)$ (resp. $H_{*}(\bar{pt})$) contains a unit $1_{pt}$ (or $1_{\bar{pt}}$) such that
$\alpha \bullet 1_{pt} = \alpha$
for all $\alpha \in F_{*}(X)$ (and similarly for $1_{\bar{pt}}$).\\
  
Let in addition  $c_{*}: F_{*}\to H_{*}$ be a natural transformation of
the associated {\em covariant} theories, with $c_{*}(1_{pt})= \bar{1}_{\bar{pt}} \in 
H_{*}(\bar{pt})$ the corresponding unit.
We are searching for a subclass of the {\em allowable} maps in $C$, containing $X\to pt$ for all objects $X$,
and a {\em partial subtheory} $F'$ of $F$ together with a {\em partial Grothendieck transformation}
\[\gamma : F'\to H\]
such that $F'(X\to pt)=F(X\to pt)$ for all objects $X$ and $\gamma_{*}=c_{*}$.\\

If one can find such a partial extension $\gamma$ of $c_{*}$, then $c_{*}$ has to {\em commute} with
suitable {\em external products}. Consider an independent square (in $C$)
\begin{displaymath} \begin{CD}
X' @>g' >> X \\
@V f' VV  @V f V \alpha V \\
Y' @> g > \beta > pt \:,
\end{CD} \end{displaymath}
with $f'$ {\em allowable} and define the  {\em external product}
\[ \times: F_{*}(Y')\times F_{*}(X) \to F_{*}(X')\]
as in \cite[p.24]{FM} by 
\[\beta \times \alpha := g^{*}(\alpha)\bullet \beta\:.\]
If in addition $f'$ is also {\em allowable} with respect to $F'$
(so that the pull-back $g^{*}$ maps $F_{*}(X)=F'_{*}(X)$
to $F'_{*}(X'\to Y')\subset F_{*}(X'\to Y')$), then one gets
\[c_{*}(\beta \times \alpha) = \gamma_{*}(\beta \times \alpha) = \gamma_{*}(g^{*}(\alpha)\bullet \beta) =\]
\[\gamma_{f'}(g^{*}(\alpha))\,\bar{\bullet}\, \gamma_{*}(\beta) = 
\bar{g}^{*}(\gamma_{f}(\alpha))\,\bar{\bullet}\, \gamma_{*}(\beta) =\]
\[ \gamma_{*}(\beta) \,\bar{\times}\, \gamma_{*}(\alpha)
= c_{*}(\beta) \,\bar{\times}\, c_{*}(\alpha) \:.\]

Next we want to {\em use} the uniqueness results of the beginning of this paper. \\

Fix a class of objects $Y$ of
$C$, called {\em orientable} (with respect to $c_{*}$), containing the final object $pt$ such that 
$F_{*}(Y)$ contains a
distinguished element $e_{Y}$ with $c_{*}(e_{Y})$ a {\em strong orientation} in $H$
(and $e_{pt}=1_{pt}\in F_{*}(pt)$). Then we call a morphism $f: X\to Y$ in $C$ {\em o-allowable}
iff $f$ is {\em allowable} (with respect to $F$) and 
the target $Y$ is {\em orientable}. We define for such an {\em o-allowable} morphism the transformation
\[\gamma_{f}: F(f:X\to Y) \to H(\bar{f}:\bar{X}\to \bar{Y})\]
as in (\ref{eq:product}) by
\begin{equation} \label{eq:defbiv}
c_{*}(\alpha\bullet e_{Y}) = \gamma_{f}(\alpha)\,\bar{\bullet}\, c_{*}(e_{Y}) \:.
\end{equation}
For $Y=pt$ we get especially $\gamma_{*}=c_{*}$, since $e_{pt}=1_{pt}$ and  
$c_{*}(1_{pt})= \bar{1}_{\bar{pt}}$. So we get a transformation $\gamma$ from $F$ restricted to these
{\em o-allowable} maps to $H$, but this need {\em not} to be a partial Grothendieck transformation.\\

Newertheless, this transformation commutes automatically with {\em pushdown}.
Consider two morphisms $g: X\to Z$ and $h: Z\to Y$ in $C$, with $g$ {\em confined} and $h$ {\em o-allowable}.
Then one gets for $f:=h\circ g : X\to Y$:
\begin{equation} \label{eq:pushdown}
\bar{g}_{*}\circ \gamma_{f} = \gamma_{h}\circ g_{*} : 
F(f: X\to Y) \to H(\bar{f}: \bar{X}\to \bar{Y}) \:,
\end{equation}
because 
\[\gamma_{h}\Bigl(g_{*}(\alpha)\Bigr)\,\bar{\bullet}\, c_{*}(e_{Y}) \stackrel{(8)}{=}
 c_{*}\Bigl(g_{*}(\alpha)\bullet e_{Y}\Bigr) \stackrel{A12}{=}
 c_{*}\Bigl(g_{*}(\alpha\bullet e_{Y})\Bigr) \stackrel{(\ast)}{=}  \]
\[ \bar{g}_{*}\Bigl(c_{*}(\alpha\bullet e_{Y})\Bigr) \stackrel{(8)}{=} 
\bar{g}_{*}\Bigl(\gamma_{f}(\alpha)\,\bar{\bullet}\, c_{*}(e_{Y})\Bigr) \stackrel{A12}{=}
\bar{g}_{*}\Bigl(\gamma_{f}(\alpha)\Bigr) \,\bar{\bullet}\, c_{*}(e_{Y})\:.\]

$ $\\
Here the equality $(\ast)$ comes from the functoriality of $c_{*}$.\\

The commutativity with {\em product} and {\em pull-back} will be build in by the following

\begin{defi} \label{defi:subgroup}
For an {\em o-allowable} morphism $f: X\to Y$ in $C$, define the subgroup
\[F'(f: X\to Y) \subset F(f: X\to Y)\]
to be the set of all $\alpha\in F(f: X\to Y)$ satisfying 
for any indepent square
\begin{displaymath} \begin{CD}
X' @>g' >> X \\
@V g^{*}(\alpha) Vf' V  @V f V \alpha V \\
Y' @> g >> Y \:,
\end{CD} \end{displaymath}
 with $f'$ also {\em o-allowable} the following two conditions:
\begin{equation} \label{eq:defproduct} 
c_{*}\Bigl(g^{*}(\alpha)\bullet \beta\Bigr) = \gamma_{f'}\Bigl(g^{*}(\alpha)\Bigr)\,\bar{\bullet}\, c_{*}(\beta) 
\end{equation}
for any $\beta\in F_{*}(Y')$, and
\begin{equation}  \label{eq:defpullback}
\gamma_{f'}\Bigl(g^{*}(\alpha)\Bigr) = \bar{g}^{*}\Bigl(\gamma_{f}(\alpha)\Bigr) \:.
\end{equation}
\end{defi}

The next lemma shows that condition (\ref{eq:defproduct}) implies also the commutativity
of $\gamma$ with {\em products} (we assume (!) that any commutative square as above with $g=id_{Y},g'=id_{X}$
is independent \cite[(B1), page 16]{FM}).

\begin{lem} \label{lem:product}
Let $f: X\to Y$ be {\em o-allowable} and assume $\alpha\in F(f: X\to Y)$ satisfies
\begin{equation} \label{eq:lemproduct} 
c_{*}(\alpha\bullet \beta) = \gamma_{f}(\alpha)\,\bar{\bullet}\, c_{*}(\beta) 
\end{equation} 
for any $\beta\in F_{*}(Y)$. Then one also has for $\beta\in F(g: Y\to Z)$, with $g$ {\em o-allowable}
the equality
\begin{equation} \label{eq:lemproduct2} 
\gamma_{g\circ f}(\alpha\bullet \beta) = \gamma_{f}(\alpha)\,\bar{\bullet}\, \gamma_{g}(\beta)\:. 
\end{equation} 
\end{lem}

This is an easy application of the definition of $\gamma$, the associativity of the bivariant product and
the assumption (\ref{eq:lemproduct}):
\[\gamma_{g\circ f}(\alpha\bullet \beta)\,\bar{\bullet}\, c_{*}(e_{Z}) \stackrel{(8)}{=}
c_{*}\Bigl((\alpha\bullet \beta)\bullet e_{Z}\Bigr) \stackrel{A1}{=}\]
\[c_{*}\Bigl(\alpha\bullet (\beta\bullet e_{Z})\Bigr) \stackrel{(12)}{=}
 \gamma_{f}(\alpha)\,\bar{\bullet}\, c_{*}(\beta\bullet e_{Z}) \stackrel{(8)}{=}\]
\[\gamma_{f}(\alpha)\,\bar{\bullet}\, \Bigl(\gamma_{g}(\beta)\,\bar{\bullet}\, c_{*}(e_{Z})\Bigr) \stackrel{A1}{=}
 \Bigl(\gamma_{f}(\alpha)\,\bar{\bullet}\,\gamma_{g}(\beta)\Bigr)\,\bar{\bullet}\, c_{*}(e_{Z}) \:.\]

Now we are ready to prove the main result of this paper
(following \cite{Y6}, and compare also with \cite[thm.A, thm.3.10]{Y5} for a similar construction
of a Grothendieck transformation to an operational bivariant theory):

\begin{thm} \label{thm:main}
$F'$ is a {\em partial subtheory} of $F$ and $\gamma: F'\to H$ is a 
{\em partial Grothendieck transformation}. Assume in addition that $c_{*}: F_{*}\to H_{*}$
commutes with {\em external products} for any independent square (in $C$)
\begin{displaymath} \begin{CD}
X' @>g' >> X \\
@V f' VV  @V f VV \\
Y' @> g >> pt \:,
\end{CD} \end{displaymath}
with $f'$ {\em o-allowable}. Then $F'_{*}=F_{*}$ and $\gamma_{*}=c_{*}$.
Moreover, $F'$ is then independent of the choice of $e_{Y}\in F_{*}(Y)$ (for $Y$ orientable).
More precisely, consider 
another partial Grothendieck transformation $\gamma'': F''\to H$, with  $F''_{*}=F_{*}$, 
$\gamma''_{*}=c_{*}$ , which is also defined on all {\em o-allowable} maps $f: X\to Y$.
Then $F''(f: X\to Y)\subset$ $F'(f: X\to Y)$ and $\gamma_{f}''=\gamma_{f}$ for all such $f$.
 \end{thm}

\begin{proof} For the first part, we only have to show that $F'$ is "stable" under
pull-back, product and push-down (then $\gamma$ is a partial Grothendieck transformation
by the results we already explained before).
\begin{enumerate}
\item {\bf pull-back}: Consider two independent squares 
\begin{displaymath} \begin{CD}
X'' @> h' >>X' @>g' >> X \\
@V f'' VV @V f' V g^{*}(\alpha) V  @V f V \alpha V \\
Y'' @> h >> Y' @> g >> Y \:,
\end{CD} \end{displaymath}
with $f,f',f''$ {\em o-allowable}. Then one gets for $\alpha\in F'(f: X\to Y)$ and for 
all $\beta\in F_{*}(Y'')$:
\[c_{*}\Bigl(h^{*}(g^{*}(\alpha))\bullet \beta\Bigr) \stackrel{A3}{=}
c_{*}\Bigl((g\circ h)^{*}(\alpha)\bullet \beta\Bigr) \stackrel{(10)}{=} \]
\[\gamma_{f''}\Bigl((g\circ h)^{*}(\alpha)\Bigr) \,\bar{\bullet}\, c_{*}(\beta) \stackrel{A3}{=}
\gamma_{f''}\Bigl(h^{*}(g^{*}(\alpha))\Bigr) \,\bar{\bullet}\, c_{*}(\beta) \:.\]
And similarly
\[\gamma_{f''}\Bigl(h^{*}(g^{*}(\alpha))\Bigr) \stackrel{A3}{=}
\gamma_{f''}\Bigl((g\circ h)^{*}(\alpha)\Bigr) \stackrel{(11)}{=}\]
\[(\overline{g\circ h})^{*}\Bigl(\gamma_{f}(\alpha)\Bigr) \stackrel{A3}{=}
\bar{h}^{*}\Bigl(\bar{g}^{*}(\gamma_{f}(\alpha))\Bigr) \stackrel{(11)}{=}
\bar{h}^{*}\Bigl(\gamma_{f'}(g^{*}(\alpha))\Bigr) \:.\]
\item {\bf product}:
Consider two independent squares 
\begin{displaymath} \begin{CD}
X' @>h'' >> X \\
@V f' VV  @V f V \alpha V \\
Y' @> h' >> Y \\
@V g' VV @V g V \beta V\\
Z' @> h >> Z \:,
\end{CD} \end{displaymath}
with the maps $f,f',g,g'$ {\em o-allowable}. Then one gets for $\alpha\in F'(f: X\to Y)$,
$\beta\in F'(f: Y\to Z)$ and for all $\delta\in F_{*}(Z')$:
\[c_{*}\Bigl(h^{*}(\alpha\bullet\beta)\bullet \delta\Bigr) \stackrel{A13}{=}
c_{*}\Bigl((h'^{*}(\alpha)\bullet h^{*}(\beta))\bullet \delta\Bigr) \stackrel{A1}{=}\]
\[c_{*}\Bigl(h'^{*}(\alpha)\bullet (h^{*}(\beta)\bullet \delta)\Bigr) \stackrel{(10)}{=}
\gamma_{f'}\Bigl(h'^{*}(\alpha)\Bigr) \,\bar{\bullet}\, c_{*}\Bigl(h^{*}(\beta)\bullet \delta\Bigr) 
 \stackrel{(10)}{=}\]
\[\gamma_{f'}\Bigl(h'^{*}(\alpha)\Bigr) \,\bar{\bullet}\, \Bigl(\gamma_{g'}(h^{*}(\beta))
\,\bar{\bullet}\, c_{*}(\delta)\Bigr) \stackrel{A1}{=}\]
\[\Bigl(\gamma_{f'}(h'^{*}(\alpha)) \,\bar{\bullet}\, \gamma_{g'}(h^{*}(\beta))\Bigr)
\,\bar{\bullet}\, c_{*}(\delta) \:.\]
But $h'^{*}(\alpha)\in F'(f': X'\to Y')$ by 1. so that we can apply lemma \ref{lem:product} and the above
equalities continue as follows:
\[\stackrel{(13)}{=} 
\gamma_{g'\circ f'}\Bigl((h'^{*}(\alpha)\bullet h^{*}(\beta)\Bigr)\,\bar{\bullet}\, c_{*}(\delta) 
\stackrel{A13}{=}
\gamma_{g'\circ f'}\Bigl(h^{*}(\alpha\bullet\beta)\Bigr)\,\bar{\bullet}\, c_{*}(\delta)\:.\]
And similarly (again using lemma \ref{lem:product})
\[\gamma_{g'\circ f'}\Bigl(h^{*}(\alpha\bullet\beta)\Bigr) \stackrel{A13}{=}
\gamma_{g'\circ f'}\Bigl(h'^{*}(\alpha)\bullet h^{*}(\beta)\Bigr) \stackrel{(13)}{=}\]
\[\gamma_{f'}(h'^{*}(\alpha)) \,\bar{\bullet}\, \gamma_{g'}(h^{*}(\beta)) \stackrel{(11)}{=}
\bar{h'}^{*}(\gamma_{f}(\alpha)) \,\bar{\bullet}\, \bar{h}^{*}(\gamma_{g}(\beta)) \stackrel{A13}{=}\]
\[\bar{h}^{*}\Bigl( \gamma_{f}(\alpha)) \,\bar{\bullet}\, \gamma_{g}(\beta) \Bigr) \stackrel{(13)}{=}
\bar{h}^{*}\Bigl( \gamma_{g\circ f}(\alpha\bullet \beta) \Bigr) \:.\]
\item {\bf push-down}:
Consider the independent squares as in 2., with $g,g'$ {\em o-allowable}.
Assume $f$ (and therefore also $f'$) is {\em confined} and fix $\alpha\in F'(g\circ f: X\to Z)$.
Then one gets for all $\beta\in F_{*}(Z')$:
\[c_{*}\Bigl(h^{*}(f_{*}(\alpha))\bullet \beta\Bigr) \stackrel{A23}{=}
c_{*}\Bigl(f'_{*}(h^{*}(\alpha))\bullet \beta\Bigr) \stackrel{A12}{=}\]
\[c_{*}\Bigl(f'_{*}(h^{*}(\alpha)\bullet \beta)\Bigr) \stackrel{(\ast)}{=}
\bar{f}'_{*}\Bigl(c_{*}(h^{*}(\alpha)\bullet \beta)\Bigr) \stackrel{(10)}{=}\]
\[\bar{f}'_{*}\Bigl(\gamma_{g'\circ f'}(h^{*}(\alpha))\,\bar{\bullet}\, c_{*}(\beta)\Bigr) \stackrel{A12}{=}
\bar{f}'_{*}\Bigl(\gamma_{g'\circ f'}(h^{*}(\alpha))\Bigr)\,\bar{\bullet}\, c_{*}(\beta) \stackrel{(9)}{=}\]
\[\gamma_{g'}\Bigl(f'_{*}(h^{*}(\alpha))\Bigr)\,\bar{\bullet}\, c_{*}(\beta) \stackrel{A23}{=}
\gamma_{g'}\Bigl(h^{*}(f_{*}(\alpha))\Bigr)\,\bar{\bullet}\, c_{*}(\beta)\:.\]
Here $(\ast)$ follows from functoriality of $c_{*}$.
Finally we have again by (\ref{eq:pushdown}):
\[\gamma_{g'}\Bigl(h^{*}(f_{*}(\alpha))\Bigr) \stackrel{A23}{=}
\gamma_{g'}\Bigl(f'_{*}(h^{*}(\alpha))\Bigr) \stackrel{(9)}{=}\]
\[\bar{f}'_{*}\Bigl(\gamma_{g'\circ f'}(h^{*}(\alpha))\Bigr) \stackrel{(11)}{=}
\bar{f}'_{*}\Bigl(\bar{h}^{*}(\gamma_{g\circ f}(\alpha))\Bigr) \stackrel{A23}{=}\]
\[\bar{h}^{*}\Bigl(\bar{f}_{*}(\gamma_{g\circ f}(\alpha))\Bigr) \stackrel{(9)}{=}
\bar{h}^{*}\Bigl(\gamma_{g}(f_{*}(\alpha))\Bigr) \:.\]
\end{enumerate}
$ $\\

Assume now that $c_{*}$ commutes with {\em external products} as stated in the theorem.
Then one gets for any independent square
\begin{displaymath} \begin{CD}
X' @>g' >> X \\
@V f' VV  @V f V \alpha V \\
Y' @> g > \beta> pt \:,
\end{CD} \end{displaymath}
with $f'$ {\em o-allowable}:
\begin{equation} \label{eq:external}
c_{*}(\beta \times \alpha) = c_{*}(\beta) \,\bar{\times}\, c_{*}(\alpha) 
\end{equation}
for $\alpha\in F_{*}(X)$ and $\beta\in F_{*}(Y')$. 
For $\beta=e_{Y'}$ one gets especially
\[c_{*}(g^{*}(\alpha)\bullet e_{Y'}) = \bar{g}^{*}(c_{*}(\alpha)) \,\bar{\bullet}\, c_{*}(e_{Y'}) \]
and this implies by definition (i.e. by equation (\ref{eq:defbiv})):
\[\gamma_{f'}(g^{*}(\alpha)) = \bar{g}^{*}(c_{*}(\alpha))\:.\]
So $\alpha$ satisfies (\ref{eq:defpullback}), and  together with
(\ref{eq:external}) this implies also the property (\ref{eq:defproduct}), i.e.
\[F'_{*}(X) = F_{*}(X) \:.\]
The last statement of the theorem follows directly from the fact that the corresponding
bivariant transformation
\[\gamma_{f}'': F''(f:X\to Y) \to H(\bar{f}:\bar{X}\to \bar{Y})\]
is uniquely defined by
\[c_{*}(\alpha\bullet e_{Y}) = \gamma_{f}''(\alpha)\,\bar{\bullet}\, c_{*}(e_{Y}) \:,\]
since $c_{*}(e_{Y})$ is a {\em strong orientation}. Especially, two such subtheories
$F',F''$ defined by different choices of the $e_{Y}\in F_{*}(Y)$ (for $Y$ orientable) have to agree.
\end{proof}

\begin{rem} The last part of theorem \ref{thm:main} shows especially, that one gets for
$c_{*}=\gamma_{*}'': F_{*}\to H_{*}$ the associated covariant functor of a
partial Grothendieck transformation $\gamma'': F\to H$ of partial (weak) bivariant theories:
\[\gamma = \gamma'': F'=F \to H\:,\]
with $F$ restricted to all {\em o-allowable} maps.
So one gets in this case nothing new (as it should be).
\end{rem}

The covariant transformation $\gamma_{*}$ in our examples comes in the following
cases from a Grothendieck transformation $\gamma: F\to H$ of bivariant theories:
\begin{itemize}
\item {\bf example \ref{ex:charclass}.1,} if all schemes are quasi-projective over a fixed
non-singular base $S$. Here $F(f)=K(f)$ is the Grothendieck group of $f$-perfect complexes and
\[\gamma:=\tau : K \to H:=A\otimes \bb{Q}\]
is the  Grothendieck transformation of \cite[thm., p.366]{Ful}.
\item {\bf example \ref{ex:charclass}.3,} if we only consider "cellular" holomorphic maps
between complex spaces, which can be embedded into smooth complex manifolds.
Here $\gamma:=c$ is the bivariant Chern transformation constructed in \cite{Br}
(and compare with \cite{Sab2} for a corresponding bivariant theory of Lagrangian cycles). 
\item  {\bf example \ref{ex:charclass}.4.} for the pl-context. Here $\gamma:=w$ is the
bivarint Stiefel-Whitney transformation of \cite[Part I, §6]{FM} (corrected in \cite{EH}).
\item {\bf example \ref{ex:cohop}.1} and {\bf example \ref{ex:LefRR}.1.}
\item Finally also the transformation $\gamma_{*}$ in the context of
{\em oriented cohomology (pre)theories} comes from a bivariant transformation
$\gamma: A\to B$ defined for a morphism $f: X\to Y$ of smooth manifolds as
\[\gamma_{f}:=\phi\,\cup td_{\phi}(TX)\cup f^{*}(td_{\phi}(TY))^{-1}\:.\]
Here $\phi=\gamma^{*}: A=A^{*}\to B=B^{*}$ is the given "nice" ring morphism.
That $\gamma$ commutes with the bivariant product follows from the projection formula
and the fact, that $\phi$ is ring morphism. For the commutativity with push-down one has in addition
to use the Riemann-Roch theorem \cite[thm.2.5.4]{Pa}.
 That $\gamma$ commutes with pull-back follows finally
from the fact, that all independent squares are transverse cartesian diagrams, e.g.
the class 
\[TX - f^{*}TY\:\in K^{0}(X)\]
behaves well under pullback in transverse cartesian diagrams
(and $td_{\phi}: K^{0}(X)\to B(X)$ commutes with pullback \cite[prop.2.2.3]{Pa}).
\end{itemize}

Before we apply theorem \ref{thm:main} to the rest of our previous examples, let us recall the main 
formal properties of bivariant theories \cite{FM}, which remain true (with slight modifications)
in the context of a {\em partial (weak) bivariant theory} $F$ on the category $C\,$:

\begin{itemize}
\item We assume that all maps $X\to pt$ are {\em allowable} so that one has an associated
{\em covariant} theory (functorial with repect to {\em confined} maps).
\item For objects $X$ of $C$ with $id: X\to X$ {\em allowable}, one has an associated
group $F^{*}(X):=F(id:X\to X)$, with a {\em cup-product} (given by the bivariant product
for the composite of this identity morphisms). This is {\em contravariant} with respect
to the class of {\em co-confined} maps, i.e. maps $f:X\to Y$ such that the square
\begin{displaymath} \begin{CD}
X @>f>> Y \\
@V id_{X} VV  @VV id_{Y} V \\
X @> f >> Y 
\end{CD} \end{displaymath}
is independent (with $id_{X},id_{Y}$ allowable). Moreover for such $X$ one also has a
{\em cap-product}
\[\cap\,: F^{*}(X)\times F_{*}(X)\to F_{*}(X)\:.\]
\item Similarly, for $f: X\to Y$ and $id_{Y}: Y\to Y$ allowable, one has a {\em right action}
\[\cap f^{*} := \bullet\,: F(f:X\to Y)\times F^{*}(Y)\to F(f:X\to Y)\:,\]
which makes $F(f:X\to Y)$ into an (unitary) right $F^{*}(Y)$-module.
\item For an indepedent square
\begin{displaymath} \begin{CD}
X' @>>> X \\
@V f' VV  @V f V \alpha V \\
Y' @> g > \beta > pt 
\end{CD} \end{displaymath}
with $f'$ {\em allowable}, one has an {\em external product}
\[ \times: F_{*}(Y')\times F_{*}(X) \to F_{*}(X')\]
defined by 
$\beta \times \alpha := g^{*}(\alpha)\bullet \beta\:$. 
\item For an independent square
\begin{displaymath} \begin{CD}
X_{pt} @>>> X \\
@VVV  @V f V \alpha V \\
pt @> i >> Y 
\end{CD} \end{displaymath}
with $f$ {\em allowable}, one has a {\em restriction to the fiber}:
\[i^{*}: F(f:X\to Y)\to F_{*}(X_{pt})\:.\]
\item Any element $\theta\in F(f:X\to Y)$ for $f$ allowable induces {\em Gysin homomorphisms}
\[\theta^{*}: F_{*}(Y)\to F_{*}(X)\:,\:\: \theta^{*}(\alpha):=\theta\bullet \alpha\:,\]
and for $f$ also confined
\[\theta_{*}: F^{*}(X)\to F^{*}(Y)\:,\:\: \theta_{*}(\beta):= f_{*}(\beta\bullet\theta) \:.\]
These Gysin homomorphisms are {\em functorial} in $\theta$ \cite[(G1), p.26]{FM}.
\end{itemize} 

If in addition $\gamma: F\to H$ is a {\em partial Grothendieck transformation}
of partial (weak) bivariant theories, then $\gamma$ induces associated transformations
$\gamma^{*}: F^{*}\to H^{*}$ and $\gamma_{*}: F_{*}\to H_{*}$, with the {\em module properties}
\[\gamma^{*}(\alpha)\cap \gamma_{*}(\beta) = \gamma_{*}(\alpha\cap \beta)\]
and
\[\gamma_{f}(\alpha)\cap \bar{f}^{*} (\gamma_{*}(\beta)) = \gamma_{f}(\alpha\cap f^{*}(\beta))\:,\]
if $\cap$ or $\cap f^{*}$ is defined. $\gamma_{*}$ commutes also with 
{\em exterior products} (if defined):
\[\gamma_{*}(\alpha\times \beta) = \gamma_{*}(\alpha)\times \gamma_{*}(\beta)\:,\]
{\em restriction to fibers} and {\em Gysin homomorphisms}.\\

\section{Examples}
 
Let us come back to the given examples of our paper. We would like to apply
theorem \ref{thm:main}. First we remark, that in all examples $\gamma_{*}$ commutes
with exterior products:
\begin{itemize}
\item For example \ref{ex:charclass}.1 this follows from \cite[Ex. 18.3.1, p.360]{Ful}.
\item For example \ref{ex:charclass}.2-3 this follows from \cite{Sch} in the case of spaces
that can be embedded into smooth manifolds. The general case follows by resolution of singularities
(compare \cite{Kw,KwY}). 
\item For example \ref{ex:charclass}.4-5 in the subanalytic context this follows from \cite{Sch}.
Here all subanalytic sets are assumed to be given in a real analytic manifold. 
The pl-context follows from the corresponding bivariant theory \cite{FM,EH}, and goes back to
\cite{HT}. 
\item For example \ref{ex:cohop}.2 this is not explicitly stated in \cite[p.190]{FL},
but follows from the skeched construction there. In the form sufficient for our applications
(i.e. one factor is a smooth manifold), it follows also from {\bf SSR 1-3} of \cite[p.188-190]{FL}.
\item For example \ref{ex:cohop}.3 this follows from \cite[prop.10.4(ii)]{Bro}.
\item Finally  example \ref{ex:LefRR}.2 follows from the corresponding property of the 
Chern class transformation $c_{*}$ already explained above, and the simple fact
that the transformation 
\begin{displaymath} \begin{CD}
K^{eq}_{c}(X\to pt) @> tr_{|X|}\circ \,|\,>> F(|X|)\otimes k 
\end{CD} \end{displaymath}
commutes with exterior products . 
\end{itemize}

The final piece of information that we need about our examples, is the fact
that for the covariant transformations $\gamma_{*}: F_{*}\to H_{*}$ given there, 
 $F_{*},H_{*}$ are just the associated covariant theories of suitable (weak) bivariant theories
 $F,H$ (together with a corresponding functor $\:\bar{ }\:$ of the underlying categories). 
For $H_{*}$ this is already the case. So we only have to deal with $F_{*}$ in the
 cases, where $\gamma_{*}$ is not already induced from a Grothendieck transformation
 $\gamma$ of bivariant theories. For all these remaining cases, the following general
 construction of (what we call) {\em simple} (weak) bivariant theories applies
 (in the context of constructible functions this notion goes back to \cite{Y2,Y5}).\\

Let $C$ be a category with classes of confined maps, independent squares
and a final object (as in \cite[Part I, §2]{FM}). We make in addition the following assumptions:

\begin{itemize}
\item {\bf (SB1)} We have a contravariant functor
\[F: C \to Rings \]
  with values in the category of rings with unit.
\item {\bf (SB2)} $F$ is also covariant functorial with respect to the confined maps
(as a functor to the category of Abelian groups).
\item {\bf (SB3)} $F$  satisfies the projection-formula  (i.e. the push-down for $f: X\to Y$ confined
is a right $F(Y)$-module operator).
\item {\bf (SB4)} $F$ has the base-change property 
\[g^{*}f_{*}=f'_{*}g'^{*} : \: F(X)\to F(Y')\]
for any independent square
\begin{displaymath} \begin{CD}
X' @>g' >> X \\
@V f' VV  @VV f V \\
Y' @> g >> Y \:,
\end{CD} \end{displaymath}
with $f,f'$ confined.
\end{itemize}  
   
The example we have already seen was the case of an 
{\em oriented cohomology (pre)theory}. And as in that case one gets a {\em simple} 
weak bivariant theory
$F$ by $F(X) := F(f:X\to Y)$, with the obvious push-down and pull-back transformations.
Finally the bivariant product
\[\bullet: F(f:X\to Y) \times F(g:Y\to Z) \to F(g\circ f: X\to Z)\]
is just given by $\alpha\bullet \beta := \alpha \cup f^{*}(\beta)$, with $\cup$ the given product
of the ring-structure.
We leave it to the reader to check that this defines a weak bivariant theory (with units)
in the sense of \cite{FM} (i.e. without property $(A123)$).
Suppose in addition:
\begin{itemize}
\item {\bf (SB5)} A commutative square is independent iff
its transpose \cite[p.17]{FM} is independent, and $F$ satisfies the two-sided
projection formula (i.e. the push-down for $f: X\to Y$ confined
is a two-sided $F(Y)$-module operator). 
\end{itemize}

Then the simple theory $F$ satisfies also
the projection formula  $(A123)$ of \cite{FM}. 

\begin{rem}
These simple (weak) bivariant theories are (skew-)commutative (and graded), if $F$ is functor to the
category of (skew-)commutative (graded) rings (and the push-down for confined maps is
degree preserving). Moreover, 
\[1_{f}:=1_{X}\in F(X)=:F(f:X\to Y)\]
is allways a canonical and strong orientation.
\end{rem}

In the very special case, that one considers only {\em trivial} independent squares:
\begin{displaymath} \begin{CD}
X @>id >> X \\
@V f VV  @VV f V \\
Y @> id >> Y 
\end{CD} \end{displaymath}
and {\em all} maps are confined, our assumptions {\bf (SB1-3)} just reduce to the properties {\bf F1-F3}
of \cite[p.28]{FL}.\\

Another more important example for us comes from the theory of constructible sheaves and functions
in the complex algebraic or (real sub-) analytic context, or in the algebraic context of seperated schemes  
of finite type over a field $k$ of characteristic zero. Here the confined maps are the proper
maps, and the independent squares are given by the cartesian diagrams. In the (sub)analytic context
we assume (for simplicity) that all spaces are of bounded dimension.\\

Then one can work in all cases with the corresponding "bounded derived category of constructible
sheaves" $D^{b}_{c}(X)$. Here we consider constructible sheaves of vector-spaces over a (suitable)
field  $R$ (in the algebraic context), with finite dimensional stalks (in the closed points).
These are stable under the usual pull-back $f^{*}$, the exact
tensor-product $\otimes_{R}$ and for proper $f$ also under push-down $Rf_{*}$.
These are related by a "projection-formula" and "proper base-change formula", which imply that
the functor $F$ given by the Grothendieck group 
\[F(X):=K_{0}(D^{b}_{c}(X))\]
with it induced transformations $f^{*},\otimes,f_{*}$ satisfies our assumptions {\bf (SB1-5)}
(with the unit given by the
class of the constant sheaf $R_{X}$). This is also the Grothendieck group of the abelian category
of constructible sheaves.
Moreover, by taking stalkwise the Euler-characteristic
($mod\:\: 2$ in the subanalytic context), one gets a natural surjective transformation 
\[\chi_{X}: K_{0}(D^{b}_{c}(X))\to CF(X)\]
onto the group of constructible functions (on the set of closed points in the algebraic context).
This induces the corresponding transformations (compare \cite[sec.2.3]{Sch3})
\[f^{*},\:\cdot,\:f_{*}\] 
on $CF$, which therefore
also satisfy our assumptions {\bf (SB1-5)}. A similar reasoning applies also to the equivariant context
studied in example \ref{ex:charclass}.5 and \ref{ex:LefRR}.2. The group of constructible
functions invariant under the real structure in example \ref{ex:charclass}.5 is stable under the
transformations   $f^{*},\cdot,f_{*}$ for equivariant maps, and in  example \ref{ex:LefRR}.2
one has similar transformations on the corresponding Grothendieck group $K^{eq}_{c}(X)$
of "equivariant constructible sheaves". These simple  bivariant theories are all commutative.\\

So the simple bivariant group $CF$ of (invariant) constructible functions can be used in 
example \ref{ex:charclass}.2-4 (or 5), and the simple bivariant (Grothendieck) group $K^{eq}_{c}$
can be used in example \ref{ex:LefRR}.2. \\

But one {\em cannot} use in example \ref{ex:charclass}.4
the simple bivariant group $CF$ of $\bb{Z}_{2}$-valued subanalytically constructible functions,
because the Stiefel-Whitney transformation $w_{*}$ is only defined for those constructible
functions, which are in addition {\em self-dual} (i.e. Euler functions \cite[def.4.2, p.823]{FuMC}).
But this self-duality condition is  only stable under proper push-down, but in general not under 
the transformations $f^{*}$ and $\cdot$. \\

This is one of the reasons for introducing the
example \ref{ex:charclass}.5 (which is better behaved). 
Newertheless, if we restrict the simple category
$CF$ to the subcategory of {\em smooth} subanalytic maps (i.e. submersions)
between real analytic {\em manifolds}, then the corresponding subcategory $CF_{eu}$ of Euler constructible
functions is stable under the bivariant product, push-down and pull-back and defines on this
restricted category a suitable bivariant theory for  example \ref{ex:charclass}.4.

\begin{rem} Of course one can also make further restrictions on the simple bivariant theory
$CF$ of constructible functions in the algebraic or complex analytic context.
The restriction to {\em smooth} holomorphic maps between {\em complex manifolds} was for example
used in \cite{Y2}.

 One gets another kind of restriction, if one makes additional
assumptions on the base-change map $g$ in the definition of independent squares:
Here one can for example assume that $g$ belongs to a class of morphisms (containing all identity
maps), which is stable under composition and base-change, e.g. the class of smooth or \'{e}tale
morphisms, open inclusions or projection of products. The case of {\em smooth} maps $g$ for independent
squares was used in the complex analytic context in \cite{Y3}.
\end{rem}

So we can apply theorem \ref{thm:main} in all the cases above to the corresponding
transformation $\gamma_{*}$, and the next natural question is, which
$\alpha \in F(f:X\to Y)$ (with $Y$ smooth) belong to the bivariant subgroup
$F'(f:X\to Y)$ constructed in theorem \ref{thm:main}. \\

For example \ref{ex:charclass}.3 one has the {\em conjecture} (\cite{Y2,Y3,Y6}):
\[F_{eu}(f:X\to Y)\subset F'(f:X\to Y)\:,\]
with $F_{eu}(f:X\to Y)$ the group of constructible functions satisfying a suitable
{\em local Euler condition} (compare \cite{Br,FM,Sab2,EY1,EY2,Zhou,Zhou2}).\\

Here we restrict ourself to
the important special case of the function $e_{f}:=1_{X}\in F(f:X\to Y)$ (or the class $e_{f}$ of the
constant sheaf, with its canonical isomorphism $\phi$, in the context of
example \ref{ex:LefRR}.2) for a {\em smooth} morphism $f$ (so that $X$ is also a manifold).\\

By remark \ref{rem:verdier} we know already the formula
\[\gamma_{f}(e_{f}) = c^{*}(T_{\bar{f}}) \,\bar{\bullet}\,  [\bar{f}] \;,\]
with $T_{\bar{f}}$ the relative tangent bundle (i.e. in our context of a smooth map
$\bar{f}$ this is just the class of the tangent bundle to the fibers in the
Grothendieck group $K_{0}(\bar{X})$ of vector bundles on $\bar{X}$) and
$[\bar{f}]\in H(\bar{f}: \bar{X}\to \bar{Y})$ the relative orientation class of the smooth
morphism $\bar{f}$. Finally 
\[c^{*}: K_{0}(\bar{X})\to H^{*}(\bar{X})\]
is the corresponding characteristic class (Chern or Stiefel-Whitney class,
depending on the example).
This characteristic class behaves well under pullback in cartesian diagrams,
since 
\[T_{\bar{f'}} = \bar{g}^{*}T_{\bar{f}}\]
for $f$ smooth. Therefore $e_{f}$ satisfies condition (\ref{eq:defpullback}) of
definition \ref{defi:subgroup}. Moreover, the other condition (\ref{eq:defproduct})
is equivalent to the {\em Verdier Riemann-Roch formula} (i.e. the commutativity of the following diagram)
\begin{displaymath} \begin{CD}
F_{*}(Y) @>\gamma_{*} >> H_{*}(\bar{Y}) \\
@V f^{*} VV  @VV c^{*}(T_{\bar{f}})\cap \bar{f}^{!} V \\
F_{*}(X) @>\gamma_{*} >> H_{*}(\bar{X}) \:,
\end{CD} \end{displaymath}
with $\bar{f}^{!}$ the corresponding smooth pull-back in homology.\\

By \cite{Sch,Sch2} and \cite[thm.2.2]{Y4}, 
this Verdier Riemann-Roch formula is true in all our cases (for $\bar{f}$ smooth).
So we get $e_{f}\in F'(f:X\to Y)$, and this defines 
in all examples a {\em canonical orientation} on the class
of smooth maps (between manifolds).

\begin{rem} This {\em Verdier Riemann-Roch formula} is also true for {\em smooth morphisms of singular
spaces} (e.g. by resolution of singularities this can be reduced to the case of manifolds
as in the proof of  \cite[thm.2.2]{Y4}). 
Therefore all $\alpha\in F(X) = CF(f:X\to Y)$ satisfy the condition (\ref{eq:defpullback}) of
definition \ref{defi:subgroup} for any independent square with a smooth base change map $g$. 
So if one defines in the examples above the independent squares as cartesian diagrams
with a smooth pullback map, then one only has to check the condition
(\ref{eq:defproduct}) of definition \ref{defi:subgroup} (compare \cite{Y2,Y3}).
\end{rem}

In a sequel to this paper we will explain a general construction of 
{\em partial bivariant theories}, which applies in particular directly to (equivariant)
Chow-groups \cite{Ful, EG1},
oriented Borel-Moore homology theories \cite{Lev}, equivariant K-theory \cite{BFQ,EG2}
or higher algebraic K-theory \cite{TT}. Moreover, we will illustrate the relation of our main theorem
\ref{thm:main} to corresponding known Riemann-Roch theorems as in \cite{BFM,BFQ,EG2,Gi,Sou}.

$ $\\
J\"{o}rg Sch\"{u}rmann\\
Westf. Wilhelms-Universit\"{a}t\\
SFB 478
"Geometrische Strukturen in der Mathematik" \\
Hittorfstr.27\\
48149 M\"{u}nster\\
Germany\\
E-mail: jschuerm@math.uni-muenster.de

\end{document}